\documentclass{amsart}

\usepackage{amssymb,amsmath}

\newtheorem{theorem}{Theorem}
\newtheorem{lemma}{Lemma}
\newcommand{\bt}{\begin{theorem}}
\newcommand{\et}{\end{theorem}}
\newcommand{\bl}{\begin{lemma}}
\newcommand{\el}{\end{lemma}}
\newcommand{\beq}{\begin{equation}}
\newcommand{\eeq}{\end{equation}}
\newcommand{\N}{\ensuremath{\mathbf N}}
\newcommand{\Q}{\ensuremath{\mathbf Q}}
\newcommand{\R}{\ensuremath{\mathbf R}}

\begin{document}
\title{Decomposing sequences into monotonic  subsequences}

\author{Melvyn B. Nathanson}
\address{Lehman College (CUNY),Bronx, New York 10468, and CUNY Graduate Center, New York, New York 10016}
\email{melvyn.nathanson@lehman.cuny.edu}

\author{Rohit Parikh}
\address{Brooklyn College (CUNY),
Brooklyn, New York 11210, and CUNY Graduate Center, New York, New York 10016}
\email{rparikh@gc.cuny.edu}

\author{Samer Salame}
\address{Department of Mathematics\\Bronx Community College (CUNY)\\Bronx, New York 10453}
\email{ssalame@gmail.com}

\thanks{The work of M.B.N. was supported in part by grants from the NSA Mathematical Sciences Program and the PSC-CUNY Research Award Program.  The work of R. P. was supported in part by a grant from the PSC-CUNY Research Award Program.}

\date{\today}

\begin{abstract}
The function $f:X \rightarrow Y$ is called $k$-monotonically increasing if there is a partition $X = \cup_{i=1}^k X_i$ such that $f|X_i : X_i \rightarrow Y$ is monotonically increasing for $i=1,\ldots,k$.   It is proved that a one-to-one function $f:\N \rightarrow \N$ is $k$-monotonically increasing if and only if every set of $k+1$ positive integers contains two integers $x,x'$ with $x < x'$ such that $f(x) \leq  f(x')$.  The function $f:X \rightarrow Y$ is called $k$-monotonic if there is a partition $X = \cup_{i=1}^k X_i$ such that $f|X_i : X_i \rightarrow Y$ is monotonically increasing or monotonically decreasing for $i=1,\ldots,k$.   It is also proved that there does not exist a $k$-monotonic function from \N\ onto \Q.
\end{abstract}

\maketitle

\section{$k$-monotonically increasing functions}
Let $X$ and $Y$ be partially ordered sets.  The function $f:X \rightarrow Y$ is 
\emph{monotonically increasing}  if $x,x' \in X$ and $x \leq x'$ implies $f(x) \leq f(x')$, and \emph{monotonically decreasing}  if $x,x' \in X$ and $x \leq x'$ implies $f(x) \geq f(x')$.   
The function $f$ is \emph{monotonic} if it is monotonically increasing or monotonically decreasing.

Let $k$ be a positive integer.  The function $f:X \rightarrow Y$ is called \emph{$k$-monotonically increasing} if $X$ contains subsets $X_1,\ldots, X_k$ such that $X = X_1 \cup \cdots \cup X_k$ and the functions $f|X_i: X_i \rightarrow Y$ are monotonically increasing for  $i=1,\ldots, k$, where $f|X_i$ is the function $f$ restricted to the subset $X_i$.

Consider the sequences 
\[
3,2,1,6,5,4,9,8,7,12,11,10\ldots
\]
and 
\[
1,3,2,6,5,4,10,9,8,7,15,14,13,12,11\ldots.
\]
Equivalently, consider the functions $g:\N \rightarrow \N$ and  $h:\N \rightarrow \N$ defined by
\[
g(n) = 
\begin{cases}
n+2 & \text {if $n \equiv 1 \pmod{3}$} \\
n & \text {if $n \equiv 2 \pmod{3}$ } \\
n-2 & \text {if $n \equiv 3 \pmod{3}$ }
\end{cases}
\]
and
\[
h(n) = \frac{m(m+1)}{2} - \ell
\]
where $m(m+1)/2$ is the unique triangular number such that 
\[
\frac{m(m-1)}{2} < n \leq \frac{m(m+1)}{2}
\]
and
\[
\ell = \frac{m(m+1)}{2}-n.
\]
The function $g$ is $3$-monotonically increasing, since $g$ is increasing on the sets $X_i = \{ n\in \N : n \equiv i \pmod{3}\}$ for $i = 1,2,3.$
On the other hand, if a function $f$ is $k$-monotonically increasing, then the pigeon-hole principle implies that  any set of $k+1$ elements of $X$ contains two elements $x,x'$ such that $x < x'$ and $f(x) \leq f(x')$.  This is a necessary condition for a function to be $k$-monotonically increasing.  Since the function $h$ is strictly decreasing on arbitrarily long intervals of consecutive integers, it follows that $h$ is not $k$-monotonic for any positive integer $k$.

We shall prove that this necessary condition is also sufficient for sequences of distinct positive integers, that is, for one-to-one functions $f: \N \rightarrow \N.$   

We define property $P(k)$ as follows:
Let $X$ and $Y$ be partially ordered sets. 
For $k \geq 2$, the function $f: X \rightarrow Y$ satisfies property $P(k)$ if, for every set $X'\subseteq X$ with $|X'|=k+1$, 
there exist $x,x' \in X'$ such that $x < x'$ and $f(x) \leq f(x')$. For totally ordered sets $X$ and $Y$, the function $f$ satisfies condition $P(k)$ if and only if $X$ does not contain an increasing sequence of $k+1$ elements $x_1 < x_2 < \cdots < x_k < x_{k+1}$ such that $f(x_1) > f(x_2) > \cdots f(x_k) > f(x_{k+1})$.

\bt    \label{mono:theorem}
Let $X$ be a finite or infinite set of positive integers, and let $f:X \rightarrow \N$ be a one-to-one function.  Let $k \geq 1$.  The function $f$ satisfies property $P(k)$ if and only if there are subsets $X_1,X_2,\ldots,X_k \subseteq X$ such that $X = X_1 \cup X_2 \cdots \cup X_k$ and $f$ restricted to $X_i$ is monotonically increasing for all $i = 1,\ldots, k$.
\et

\begin{proof}
By induction on $k \geq 1$.  The theorem holds for $k=1$, since this is simply the definition of monotonicity.

Assume the result holds for functions that satisfy condition $P(k-1)$ for some  $k \geq 2$.  Let $f$ be a function that satisfies $P(k)$.
We construct the set $X_k = \{a_n\}_{n=1}^{\infty}$ as follows:  Let $a_1 = 1$.  Given $a_n$, let $a_{n+1}$ be the smallest integer such that $f(a_{n+1}) > f(a_n)$.  Then
\[
a_1 < a_2 < \cdots < a_n < a_{n+1} < \cdots
\]
and
\[
f(a_1) < f(a_2) < \cdots < f(a_n) < f(a_{n+1}) < \cdots.
\]
Thus, $f$ is strictly increasing on the set $X_k$.  
Note that $X_k$ is finite if and only if $X$ is finite.

We shall prove that the set $B = X \setminus X_k$ satisfies condition $P(k-1)$.  If not, then there exists an increasing sequence  $b_1 < b_2 < \cdots < b_k$ of elements of $B$ such that $f(b_1) > f(b_2) > \cdots > f(b_k)$.  Since $1 = a_1 < b_1,$ there is a unique integer $n$ such that  $a_n \in X_k$ and $a_n < b_1 < a_{n+1}$.  Since the function $f$ is one-to-one and $a_{n+1}$ is the smallest integer such that $ f(a_n) < f(a_{n+1}) $, it follows that $f(a_n) > f(b_1).$  Thus,
\[
a_n < b_1 < b_2 < \cdots < b_k
\]
and
\[
f(a_n) > f(b_1) > f(b_2) > \cdots > f(b_k),
\]
which is impossible since $f$ satisfies property $P(k)$ on $X$.
Therefore, $f$ satisfies property $P(k-1)$ on $B$, and so there are subsets $X_1,X_2,\ldots,X_{k-1} \subseteq B$ such that $B = X_1 \cup X_2 \cdots \cup X_{k-1}$ and $f$ restricted to $X_i$ is increasing for $i = 1,\ldots, k-1$.  This completes the proof.
\end{proof}

\section{$k$-monotonic functions}
Let $X$ and $Y$ be partially ordered sets.  Let $k$ be a positive integer.  The function $f:X \rightarrow Y$ will be called \emph{$k$-monotonic} if $X$ contains subsets $X_1,\ldots, X_k$ such that $X = X_1 \cup \cdots \cup X_k$ and the functions $|X_i: X_i \rightarrow Y$ are monotonic (that is, monotonically increasing or monotonically decreasing) for all $i=1,\ldots, k$.  We shall prove that there does not exist a $k$-monotonic function from \N\ onto \Q.

\bt      \label{mono:theorem:NQ}
If $f:\N \rightarrow \Q$ is onto, then $f$ is not $k$-monotonic for any $k$.  
\et

\begin{proof}
Let $A_1, \ldots, A_r, B_1, \ldots, B_s$ be sets of positive integers, and let
\[
\N = A_1 \cup \cdots \cup A_r \cup B_1\cup \cdots \cup B_s.
\]
Let $f:\N \rightarrow \Q$ be a function such that $f|A_i$ is monotonically decreasing for $i=1,\ldots, r$ and  $f|B_j$ is monotonically increasing for $j=1,\ldots, s$.  Let
\[
a^{\ast} = \sup\left( \bigcup_{i=1}^r f(A_i) \right) = \sup\{f(\min(A_i)):i=1,\ldots,r  \}.
\]
Let 
\[
J_1 = \{j\in \{1,\ldots,s\} : \sup(f(B_j)) < \infty\}
\]
and
\[
b^{\ast} = \sup\left( \bigcup_{j\in J_1} f(B_j) \right).
\]  
If
\[
m^{\ast} = \max(a^{\ast},b^{\ast})
\]
then 
\[
\bigcup_{i=1}^r f(A_i) \cup \bigcup_{j \in J_1} f(B_j)  \subseteq (-\infty,m^{\ast}].
\]
Let
\[
J_2 =  \{1,\ldots,s\}\setminus J_1 = \{j \in \{1,\ldots,s\} : \sup(f(B_j)) = \infty\}.
\]
For each $j \in J_2$ there is an integer $b_j \in B_j$ such that $f(b_j) \geq m^{\ast}+1$.  Then
\[
\bigcup_{j\in J_2} \{f(b) : b\in B_j \text{ and } b \geq b_j\} \subseteq [m^{\ast}+1,\infty) .
\]
The set
\[
B^{\ast} = \bigcup_{j\in J_2} \{  b\in B_j : b < b_j\}
\]
is finite, and so the set
\[
f(\N) \cap (m^{\ast},m^{\ast}+1) \cap \Q = f(B^{\ast}) \cap (m^{\ast},m^{\ast}+1) \cap \Q
\]
is also finite.  Since nonempty open interval $(m^{\ast},m^{\ast}+1)$ contains infinitely many rational numbers, it follows that the function $f:\N\rightarrow \Q$ is not surjective.  Thus, no surjection from \N\ onto \Q\ is $k$-montonic.

We can also prove Theorem~\ref{mono:theorem:NQ} as follows:  If $\N = A_1 \cup \cdots \cup A_r \cup B_1\cup \cdots \cup B_s$ and if $f$ is monotonic on the sets $A_i$ and $B_j$, then the sets $f(A_i)$ and $f(B_j)$ each have exactly one limit point in  $\R \cup \{\pm \infty\}$.  If $q$ is a rational number that is not one of these limit points, then there is an open neighborhood $U$ of $q$ such that $U \cap f(\N) \subseteq \{q\}$, and so every rational number $q'\in U$ with $q' \neq q$ is not in the image of $f$.  This completes the proof.
\end{proof}

We define an \emph{interval of rational numbers} to be a set of the form $\mathcal{I} \cap \Q$, where $\mathcal{I}$ is an open, closed, or half-open, half-closed interval with nonempty interior.
An argument similar to that of Theorem~\ref{mono:theorem:NQ} proves the following result.

\bt
If $f$ is a function from \N\ onto an interval of rational numbers, then $f$ is not $k$-monotonic for any $k$.  
\et

It is an open problem to determine if there a simple criterion analogous to that of Theorem~\ref{mono:theorem} to determine if a one-to-one function $f: \N \rightarrow \N$ is $k$-monotonic.

\end{document}